\documentstyle[11pt]{article}
\newtheorem{theorem}{Theorem}[section]
\newtheorem{lemma}[theorem]{Lemma}

\begin{document}
%
%

\long\def\ig#1{\relax}
\ig{Thanks to Roberto Minio for this def'n.  Compare the def'n of
\comment in AMSTeX.}

\newcount \coefa
\newcount \coefb
\newcount \coefc
\newcount\tempcounta
\newcount\tempcountb
\newcount\tempcountc
\newcount\tempcountd
\newcount\xext
\newcount\yext
\newcount\xoff
\newcount\yoff
\newcount\gap%
\newcount\arrowtypea
\newcount\arrowtypeb
\newcount\arrowtypec
\newcount\arrowtyped
\newcount\arrowtypee
\newcount\height
\newcount\width
\newcount\xpos
\newcount\ypos
\newcount\run
\newcount\rise
\newcount\arrowlength
\newcount\halflength
\newcount\arrowtype
\newdimen\tempdimen
\newdimen\xlen
\newdimen\ylen
\newsavebox{\tempboxa}%
\newsavebox{\tempboxb}%
\newsavebox{\tempboxc}%

\makeatletter
\setlength{\unitlength}{.01em}%
\def\settypes(#1,#2,#3){\arrowtypea#1 \arrowtypeb#2 \arrowtypec#3}
\def\settoheight#1#2{\setbox\@tempboxa\hbox{#2}#1\ht\@tempboxa\relax}%
\def\settodepth#1#2{\setbox\@tempboxa\hbox{#2}#1\dp\@tempboxa\relax}%
\def\settokens[#1`#2`#3`#4]{%
     \def\tokena{#1}\def\tokenb{#2}\def\tokenc{#3}\def\tokend{#4}}
\def\setsqparms[#1`#2`#3`#4;#5`#6]{%
\arrowtypea #1
\arrowtypeb #2
\arrowtypec #3
\arrowtyped #4
\width #5
\height #6
}
\def\setpos(#1,#2){\xpos=#1 \ypos#2}

\def\bfig{\begin{picture}(\xext,\yext)(\xoff,\yoff)}
\def\efig{\end{picture}}

\def\putbox(#1,#2)#3{\put(#1,#2){\makebox(0,0){$#3$}}}

\def\settriparms[#1`#2`#3;#4]{\settripairparms[#1`#2`#3`1`1;#4]}%

\def\settripairparms[#1`#2`#3`#4`#5;#6]{%
\arrowtypea #1
\arrowtypeb #2
\arrowtypec #3
\arrowtyped #4
\arrowtypee #5
\width #6
\height #6
}

\def\resetparms{\settripairparms[1`1`1`1`1;500]\width 500}

\resetparms

\def\mvector(#1,#2)#3{
\put(0,0){\vector(#1,#2){#3}}%
\put(0,0){\vector(#1,#2){30}}%
}
\def\evector(#1,#2)#3{{
\arrowlength #3
\put(0,0){\vector(#1,#2){\arrowlength}}%
\advance \arrowlength by-30
\put(0,0){\vector(#1,#2){\arrowlength}}%
}}

\def\horsize#1#2{%
\settowidth{\tempdimen}{$#2$}%
#1=\tempdimen
\divide #1 by\unitlength
}

\def\vertsize#1#2{%
\settoheight{\tempdimen}{$#2$}%
#1=\tempdimen
\settodepth{\tempdimen}{$#2$}%
\advance #1 by\tempdimen
\divide #1 by\unitlength
}

\def\vertadjust[#1`#2`#3]{%
\vertsize{\tempcounta}{#1}%
\vertsize{\tempcountb}{#2}%
\ifnum \tempcounta<\tempcountb \tempcounta=\tempcountb \fi
\divide\tempcounta by2
\vertsize{\tempcountb}{#3}%
\ifnum \tempcountb>0 \advance \tempcountb by20 \fi
\ifnum \tempcounta<\tempcountb \tempcounta=\tempcountb \fi
}

\def\horadjust[#1`#2`#3]{%
\horsize{\tempcounta}{#1}%
\horsize{\tempcountb}{#2}%
\ifnum \tempcounta<\tempcountb \tempcounta=\tempcountb \fi
\divide\tempcounta by20
\horsize{\tempcountb}{#3}%
\ifnum \tempcountb>0 \advance \tempcountb by60 \fi
\ifnum \tempcounta<\tempcountb \tempcounta=\tempcountb \fi
}

\ig{ In this procedure, #1 is the paramater that sticks out all the way,
#2 sticks out the least and #3 is a label sticking out half way.  #4 is
the amount of the offset.}

\def\sladjust[#1`#2`#3]#4{%
\tempcountc=#4
\horsize{\tempcounta}{#1}%
\divide \tempcounta by2
\horsize{\tempcountb}{#2}%
\divide \tempcountb by2
\advance \tempcountb by-\tempcountc
\ifnum \tempcounta<\tempcountb \tempcounta=\tempcountb\fi
\divide \tempcountc by2
\horsize{\tempcountb}{#3}%
\advance \tempcountb by-\tempcountc
\ifnum \tempcountb>0 \advance \tempcountb by80\fi
\ifnum \tempcounta<\tempcountb \tempcounta=\tempcountb\fi
\advance\tempcounta by20
}

\def\putvector(#1,#2)(#3,#4)#5#6{{%
\xpos=#1
\ypos=#2
\run=#3
\rise=#4
\arrowlength=#5
\arrowtype=#6
\ifnum \arrowtype<0
    \ifnum \run=0
        \advance \ypos by-\arrowlength
    \else
        \tempcounta \arrowlength
        \multiply \tempcounta by\rise
        \divide \tempcounta by\run
        \ifnum\run>0
            \advance \xpos by\arrowlength
            \advance \ypos by\tempcounta
        \else
            \advance \xpos by-\arrowlength
            \advance \ypos by-\tempcounta
        \fi
    \fi
    \multiply \arrowtype by-1
    \multiply \rise by-1
    \multiply \run by-1
\fi
\ifnum \arrowtype=1
    \put(\xpos,\ypos){\vector(\run,\rise){\arrowlength}}%
\else\ifnum \arrowtype=2
    \put(\xpos,\ypos){\mvector(\run,\rise)\arrowlength}%
\else\ifnum\arrowtype=3
    \put(\xpos,\ypos){\evector(\run,\rise){\arrowlength}}%
\fi\fi\fi
}}

\def\putsplitvector(#1,#2)#3#4{
\xpos #1
\ypos #2
\arrowtype #4
\halflength #3
\arrowlength #3
\gap 140
\advance \halflength by-\gap
\divide \halflength by2
\ifnum \arrowtype=1
    \put(\xpos,\ypos){\line(0,-1){\halflength}}%
    \advance\ypos by-\halflength
    \advance\ypos by-\gap
    \put(\xpos,\ypos){\vector(0,-1){\halflength}}%
\else\ifnum \arrowtype=2
    \put(\xpos,\ypos){\line(0,-1)\halflength}%
    \put(\xpos,\ypos){\vector(0,-1)3}%
    \advance\ypos by-\halflength
    \advance\ypos by-\gap
    \put(\xpos,\ypos){\vector(0,-1){\halflength}}%
\else\ifnum\arrowtype=3
    \put(\xpos,\ypos){\line(0,-1)\halflength}%
    \advance\ypos by-\halflength
    \advance\ypos by-\gap
    \put(\xpos,\ypos){\evector(0,-1){\halflength}}%
\else\ifnum \arrowtype=-1
    \advance \ypos by-\arrowlength
    \put(\xpos,\ypos){\line(0,1){\halflength}}%
    \advance\ypos by\halflength
    \advance\ypos by\gap
    \put(\xpos,\ypos){\vector(0,1){\halflength}}%
\else\ifnum \arrowtype=-2
    \advance \ypos by-\arrowlength
    \put(\xpos,\ypos){\line(0,1)\halflength}%
    \put(\xpos,\ypos){\vector(0,1)3}%
    \advance\ypos by\halflength
    \advance\ypos by\gap
    \put(\xpos,\ypos){\vector(0,1){\halflength}}%
\else\ifnum\arrowtype=-3
    \advance \ypos by-\arrowlength
    \put(\xpos,\ypos){\line(0,1)\halflength}%
    \advance\ypos by\halflength
    \advance\ypos by\gap
    \put(\xpos,\ypos){\evector(0,1){\halflength}}%
\fi\fi\fi\fi\fi\fi
}

\def\putmorphism(#1)(#2,#3)[#4`#5`#6]#7#8#9{{%
\run #2
\rise #3
\ifnum\rise=0
  \puthmorphism(#1)[#4`#5`#6]{#7}{#8}{#9}%
\else\ifnum\run=0
  \putvmorphism(#1)[#4`#5`#6]{#7}{#8}{#9}%
\else
\setpos(#1)%
\arrowlength #7
\arrowtype #8
\ifnum\run=0
\else\ifnum\rise=0
\else
\ifnum\run>0
    \coefa=1
\else
   \coefa=-1
\fi
\ifnum\arrowtype>0
   \coefb=0
   \coefc=-1
\else
   \coefb=\coefa
   \coefc=1
   \arrowtype=-\arrowtype
\fi
\width=2
\multiply \width by\run
\divide \width by\rise
\ifnum \width<0  \width=-\width\fi
\advance\width by60
\if l#9 \width=-\width\fi
\putbox(\xpos,\ypos){#4}
{\multiply \coefa by\arrowlength
\advance\xpos by\coefa
\multiply \coefa by\rise
\divide \coefa by\run
\advance \ypos by\coefa
\putbox(\xpos,\ypos){#5} }%
{\multiply \coefa by\arrowlength
\divide \coefa by2
\advance \xpos by\coefa
\advance \xpos by\width
\multiply \coefa by\rise
\divide \coefa by\run
\advance \ypos by\coefa
\if l#9%
   \put(\xpos,\ypos){\makebox(0,0)[r]{$#6$}}%
\else\if r#9%
   \put(\xpos,\ypos){\makebox(0,0)[l]{$#6$}}%
\fi\fi }%
{\multiply \rise by-\coefc
\multiply \run by-\coefc
\multiply \coefb by\arrowlength
\advance \xpos by\coefb
\multiply \coefb by\rise
\divide \coefb by\run
\advance \ypos by\coefb
\multiply \coefc by70
\advance \ypos by\coefc
\multiply \coefc by\run
\divide \coefc by\rise
\advance \xpos by\coefc
\multiply \coefa by140
\multiply \coefa by\run
\divide \coefa by\rise
\advance \arrowlength by\coefa
\ifnum \arrowtype=1
   \put(\xpos,\ypos){\vector(\run,\rise){\arrowlength}}%
\else\ifnum\arrowtype=2
   \put(\xpos,\ypos){\mvector(\run,\rise){\arrowlength}}%
\else\ifnum\arrowtype=3
   \put(\xpos,\ypos){\evector(\run,\rise){\arrowlength}}%
\fi\fi\fi}\fi\fi\fi\fi}}

\def\puthmorphism(#1,#2)[#3`#4`#5]#6#7#8{{%
\xpos #1
\ypos #2
\width #6
\arrowlength #6
\putbox(\xpos,\ypos){#3\vphantom{#4}}%
{\advance \xpos by\arrowlength
\putbox(\xpos,\ypos){\vphantom{#3}#4}}%
\horsize{\tempcounta}{#3}%
\horsize{\tempcountb}{#4}%
\divide \tempcounta by2
\divide \tempcountb by2
\advance \tempcounta by30
\advance \tempcountb by30
\advance \xpos by\tempcounta
\advance \arrowlength by-\tempcounta
\advance \arrowlength by-\tempcountb
\putvector(\xpos,\ypos)(1,0){\arrowlength}{#7}%
\divide \arrowlength by2
\advance \xpos by\arrowlength
\vertsize{\tempcounta}{#5}%
\divide\tempcounta by2
\advance \tempcounta by20
\if a#8 %
   \advance \ypos by\tempcounta
   \putbox(\xpos,\ypos){#5}%
\else
   \advance \ypos by-\tempcounta
   \putbox(\xpos,\ypos){#5}%
\fi}}

\def\putvmorphism(#1,#2)[#3`#4`#5]#6#7#8{{%
\xpos #1
\ypos #2
\arrowlength #6
\arrowtype #7
\settowidth{\xlen}{$#5$}%
\putbox(\xpos,\ypos){#3}%
{\advance \ypos by-\arrowlength
\putbox(\xpos,\ypos){#4}}%
{\advance\arrowlength by-140
\advance \ypos by-70
\ifdim\xlen>0pt
   \if m#8%
      \putsplitvector(\xpos,\ypos){\arrowlength}{\arrowtype}%
   \else
      \putvector(\xpos,\ypos)(0,-1){\arrowlength}{\arrowtype}%
   \fi
\else
   \putvector(\xpos,\ypos)(0,-1){\arrowlength}{\arrowtype}%
\fi}%
\ifdim\xlen>0pt
   \divide \arrowlength by2
   \advance\ypos by-\arrowlength
   \if l#8%
      \advance \xpos by-40
      \put(\xpos,\ypos){\makebox(0,0)[r]{$#5$}}%
   \else\if r#8%
      \advance \xpos by40
      \put(\xpos,\ypos){\makebox(0,0)[l]{$#5$}}%
   \else
      \putbox(\xpos,\ypos){#5}%
   \fi\fi
\fi
}}

\def\topadjust[#1`#2`#3]{%
\yoff=10
\vertadjust[#1`#2`{#3}]%
\advance \yext by\tempcounta
\advance \yext by 10
}
\def\botadjust[#1`#2`#3]{%
\vertadjust[#1`#2`{#3}]%
\advance \yext by\tempcounta
\advance \yoff by-\tempcounta
}
\def\leftadjust[#1`#2`#3]{%
\xoff=0
\horadjust[#1`#2`{#3}]%
\advance \xext by\tempcounta
\advance \xoff by-\tempcounta
}
\def\rightadjust[#1`#2`#3]{%
\horadjust[#1`#2`{#3}]%
\advance \xext by\tempcounta
}
\def\rightsladjust[#1`#2`#3]{%
\sladjust[#1`#2`{#3}]{\width}%
\advance \xext by\tempcounta
}
\def\leftsladjust[#1`#2`#3]{%
\xoff=0
\sladjust[#1`#2`{#3}]{\width}%
\advance \xext by\tempcounta
\advance \xoff by-\tempcounta
}
\def\adjust[#1`#2;#3`#4;#5`#6;#7`#8]{%
\topadjust[#1``{#2}]
\leftadjust[#3``{#4}]
\rightadjust[#5``{#6}]
\botadjust[#7``{#8}]}

\def\putsquarep<#1>(#2)[#3;#4`#5`#6`#7]{{%
\setsqparms[#1]%
\setpos(#2)%
\settokens[#3]%
\puthmorphism(\xpos,\ypos)[\tokenc`\tokend`{#7}]{\width}{\arrowtyped}b%
\advance\ypos by \height
\puthmorphism(\xpos,\ypos)[\tokena`\tokenb`{#4}]{\width}{\arrowtypea}a%
\putvmorphism(\xpos,\ypos)[``{#5}]{\height}{\arrowtypeb}l%
\advance\xpos by \width
\putvmorphism(\xpos,\ypos)[``{#6}]{\height}{\arrowtypec}r%
}}

\def\putsquare{\@ifnextchar <{\putsquarep}{\putsquarep%
   <\arrowtypea`\arrowtypeb`\arrowtypec`\arrowtyped;\width`\height>}}
\def\square{\@ifnextchar< {\squarep}{\squarep
   <\arrowtypea`\arrowtypeb`\arrowtypec`\arrowtyped;\width`\height>}}
\def\squarep<#1>[#2`#3`#4`#5;#6`#7`#8`#9]{{
\setsqparms[#1]
\xext=\width                                          
\yext=\height                                         
\topadjust[#2`#3`{#6}]
\botadjust[#4`#5`{#9}]
\leftadjust[#2`#4`{#7}]
\rightadjust[#3`#5`{#8}]
\begin{picture}(\xext,\yext)(\xoff,\yoff)
\putsquarep<\arrowtypea`\arrowtypeb`\arrowtypec`\arrowtyped;\width`\height>%
(0,0)[#2`#3`#4`#5;#6`#7`#8`{#9}]%
\end{picture}%
}}

\def\putptrianglep<#1>(#2,#3)[#4`#5`#6;#7`#8`#9]{{%
\settriparms[#1]%
\xpos=#2 \ypos=#3
\advance\ypos by \height
\puthmorphism(\xpos,\ypos)[#4`#5`{#7}]{\height}{\arrowtypea}a%
\putvmorphism(\xpos,\ypos)[`#6`{#8}]{\height}{\arrowtypeb}l%
\advance\xpos by\height
\putmorphism(\xpos,\ypos)(-1,-1)[``{#9}]{\height}{\arrowtypec}r%
}}

\def\putptriangle{\@ifnextchar <{\putptrianglep}{\putptrianglep
   <\arrowtypea`\arrowtypeb`\arrowtypec;\height>}}
\def\ptriangle{\@ifnextchar <{\ptrianglep}{\ptrianglep
   <\arrowtypea`\arrowtypeb`\arrowtypec;\height>}}

\def\ptrianglep<#1>[#2`#3`#4;#5`#6`#7]{{
\settriparms[#1]%
\width=\height                         
\xext=\width                           
\yext=\width                           
\topadjust[#2`#3`{#5}]
\botadjust[#3``]
\leftadjust[#2`#4`{#6}]
\rightsladjust[#3`#4`{#7}]
\begin{picture}(\xext,\yext)(\xoff,\yoff)
\putptrianglep<\arrowtypea`\arrowtypeb`\arrowtypec;\height>%
(0,0)[#2`#3`#4;#5`#6`{#7}]%
\end{picture}%
}}

\def\putqtrianglep<#1>(#2,#3)[#4`#5`#6;#7`#8`#9]{{%
\settriparms[#1]%
\xpos=#2 \ypos=#3
\advance\ypos by\height
\puthmorphism(\xpos,\ypos)[#4`#5`{#7}]{\height}{\arrowtypea}a%
\putmorphism(\xpos,\ypos)(1,-1)[``{#8}]{\height}{\arrowtypeb}l%
\advance\xpos by\height
\putvmorphism(\xpos,\ypos)[`#6`{#9}]{\height}{\arrowtypec}r%
}}

\def\putqtriangle{\@ifnextchar <{\putqtrianglep}{\putqtrianglep
   <\arrowtypea`\arrowtypeb`\arrowtypec;\height>}}
\def\qtriangle{\@ifnextchar <{\qtrianglep}{\qtrianglep
   <\arrowtypea`\arrowtypeb`\arrowtypec;\height>}}

\def\qtrianglep<#1>[#2`#3`#4;#5`#6`#7]{{
\settriparms[#1]
\width=\height                         
\xext=\width                           
\yext=\height                          
\topadjust[#2`#3`{#5}]
\botadjust[#4``]
\leftsladjust[#2`#4`{#6}]
\rightadjust[#3`#4`{#7}]
\begin{picture}(\xext,\yext)(\xoff,\yoff)
\putqtrianglep<\arrowtypea`\arrowtypeb`\arrowtypec;\height>%
(0,0)[#2`#3`#4;#5`#6`{#7}]%
\end{picture}%
}}

\def\putdtrianglep<#1>(#2,#3)[#4`#5`#6;#7`#8`#9]{{%
\settriparms[#1]%
\xpos=#2 \ypos=#3
\puthmorphism(\xpos,\ypos)[#5`#6`{#9}]{\height}{\arrowtypec}b%
\advance\xpos by \height \advance\ypos by\height
\putmorphism(\xpos,\ypos)(-1,-1)[``{#7}]{\height}{\arrowtypea}l%
\putvmorphism(\xpos,\ypos)[#4``{#8}]{\height}{\arrowtypeb}r%
}}

\def\putdtriangle{\@ifnextchar <{\putdtrianglep}{\putdtrianglep
   <\arrowtypea`\arrowtypeb`\arrowtypec;\height>}}
\def\dtriangle{\@ifnextchar <{\dtrianglep}{\dtrianglep
   <\arrowtypea`\arrowtypeb`\arrowtypec;\height>}}

\def\dtrianglep<#1>[#2`#3`#4;#5`#6`#7]{{
\settriparms[#1]
\width=\height                         
\xext=\width                           
\yext=\height                          
\topadjust[#2``]
\botadjust[#3`#4`{#7}]
\leftsladjust[#3`#2`{#5}]
\rightadjust[#2`#4`{#6}]
\begin{picture}(\xext,\yext)(\xoff,\yoff)
\putdtrianglep<\arrowtypea`\arrowtypeb`\arrowtypec;\height>%
(0,0)[#2`#3`#4;#5`#6`{#7}]%
\end{picture}%
}}

\def\putbtrianglep<#1>(#2,#3)[#4`#5`#6;#7`#8`#9]{{%
\settriparms[#1]%
\xpos=#2 \ypos=#3
\puthmorphism(\xpos,\ypos)[#5`#6`{#9}]{\height}{\arrowtypec}b%
\advance\ypos by\height
\putmorphism(\xpos,\ypos)(1,-1)[``{#8}]{\height}{\arrowtypeb}r%
\putvmorphism(\xpos,\ypos)[#4``{#7}]{\height}{\arrowtypea}l%
}}

\def\putbtriangle{\@ifnextchar <{\putbtrianglep}{\putbtrianglep
   <\arrowtypea`\arrowtypeb`\arrowtypec;\height>}}
\def\btriangle{\@ifnextchar <{\btrianglep}{\btrianglep
   <\arrowtypea`\arrowtypeb`\arrowtypec;\height>}}

\def\btrianglep<#1>[#2`#3`#4;#5`#6`#7]{{
\settriparms[#1]
\width=\height                         
\xext=\width                           
\yext=\height                          
\topadjust[#2``]
\botadjust[#3`#4`{#7}]
\leftadjust[#2`#3`{#5}]
\rightsladjust[#4`#2`{#6}]
\begin{picture}(\xext,\yext)(\xoff,\yoff)
\putbtrianglep<\arrowtypea`\arrowtypeb`\arrowtypec;\height>%
(0,0)[#2`#3`#4;#5`#6`{#7}]%
\end{picture}%
}}

\def\putAtrianglep<#1>(#2,#3)[#4`#5`#6;#7`#8`#9]{{%
\settriparms[#1]%
\xpos=#2 \ypos=#3
{\multiply \height by2
\puthmorphism(\xpos,\ypos)[#5`#6`{#9}]{\height}{\arrowtypec}b}%
\advance\xpos by\height \advance\ypos by\height
\putmorphism(\xpos,\ypos)(-1,-1)[#4``{#7}]{\height}{\arrowtypea}l%
\putmorphism(\xpos,\ypos)(1,-1)[``{#8}]{\height}{\arrowtypeb}r%
}}

\def\putAtriangle{\@ifnextchar <{\putAtrianglep}{\putAtrianglep
   <\arrowtypea`\arrowtypeb`\arrowtypec;\height>}}
\def\Atriangle{\@ifnextchar <{\Atrianglep}{\Atrianglep
   <\arrowtypea`\arrowtypeb`\arrowtypec;\height>}}

\def\Atrianglep<#1>[#2`#3`#4;#5`#6`#7]{{
\settriparms[#1]
\width=\height                         
\xext=\width                           
\yext=\height                          
\topadjust[#2``]
\botadjust[#3`#4`{#7}]
\multiply \xext by2 
\leftsladjust[#3`#2`{#5}]
\rightsladjust[#4`#2`{#6}]
\begin{picture}(\xext,\yext)(\xoff,\yoff)%
\putAtrianglep<\arrowtypea`\arrowtypeb`\arrowtypec;\height>%
(0,0)[#2`#3`#4;#5`#6`{#7}]%
\end{picture}%
}}

\def\putAtrianglepairp<#1>(#2)[#3;#4`#5`#6`#7`#8]{{
\settripairparms[#1]%
\setpos(#2)%
\settokens[#3]%
\puthmorphism(\xpos,\ypos)[\tokenb`\tokenc`{#7}]{\height}{\arrowtyped}b%
\advance\xpos by\height
\advance\ypos by\height
\putmorphism(\xpos,\ypos)(-1,-1)[\tokena``{#4}]{\height}{\arrowtypea}l%
\putvmorphism(\xpos,\ypos)[``{#5}]{\height}{\arrowtypeb}m%
\putmorphism(\xpos,\ypos)(1,-1)[``{#6}]{\height}{\arrowtypec}r%
}}

\def\putAtrianglepair{\@ifnextchar <{\putAtrianglepairp}{\putAtrianglepairp%
   <\arrowtypea`\arrowtypeb`\arrowtypec`\arrowtyped`\arrowtypee;\height>}}
\def\Atrianglepair{\@ifnextchar <{\Atrianglepairp}{\Atrianglepairp%
   <\arrowtypea`\arrowtypeb`\arrowtypec`\arrowtyped`\arrowtypee;\height>}}

\def\Atrianglepairp<#1>[#2;#3`#4`#5`#6`#7]{{%
\settripairparms[#1]%
\settokens[#2]%
\width=\height
\xext=\width
\yext=\height
\topadjust[\tokena``]%
\vertadjust[\tokenb`\tokenc`{#6}]
\tempcountd=\tempcounta                       
\vertadjust[\tokenc`\tokend`{#7}]
\ifnum\tempcounta<\tempcountd                 
\tempcounta=\tempcountd\fi                    
\advance \yext by\tempcounta                  
\advance \yoff by-\tempcounta                 %
\multiply \xext by2 
\leftsladjust[\tokenb`\tokena`{#3}]
\rightsladjust[\tokend`\tokena`{#5}]%
\begin{picture}(\xext,\yext)(\xoff,\yoff)%
\putAtrianglepairp
<\arrowtypea`\arrowtypeb`\arrowtypec`\arrowtyped`\arrowtypee;\height>%
(0,0)[#2;#3`#4`#5`#6`{#7}]%
\end{picture}%
}}

\def\putVtrianglep<#1>(#2,#3)[#4`#5`#6;#7`#8`#9]{{%
\settriparms[#1]%
\xpos=#2 \ypos=#3
\advance\ypos by\height
{\multiply\height by2
\puthmorphism(\xpos,\ypos)[#4`#5`{#7}]{\height}{\arrowtypea}a}%
\putmorphism(\xpos,\ypos)(1,-1)[`#6`{#8}]{\height}{\arrowtypeb}l%
\advance\xpos by\height
\advance\xpos by\height
\putmorphism(\xpos,\ypos)(-1,-1)[``{#9}]{\height}{\arrowtypec}r%
}}

\def\putVtriangle{\@ifnextchar <{\putVtrianglep}{\putVtrianglep
   <\arrowtypea`\arrowtypeb`\arrowtypec;\height>}}
\def\Vtriangle{\@ifnextchar <{\Vtrianglep}{\Vtrianglep
   <\arrowtypea`\arrowtypeb`\arrowtypec;\height>}}

\def\Vtrianglep<#1>[#2`#3`#4;#5`#6`#7]{{
\settriparms[#1]
\width=\height                         
\xext=\width                           
\yext=\height                          
\topadjust[#2`#3`{#5}]
\botadjust[#4``]
\multiply \xext by2 
\leftsladjust[#2`#3`{#6}]
\rightsladjust[#3`#4`{#7}]
\begin{picture}(\xext,\yext)(\xoff,\yoff)%
\putVtrianglep<\arrowtypea`\arrowtypeb`\arrowtypec;\height>%
(0,0)[#2`#3`#4;#5`#6`{#7}]%
\end{picture}%
}}

\def\putVtrianglepairp<#1>(#2)[#3;#4`#5`#6`#7`#8]{{
\settripairparms[#1]%
\setpos(#2)%
\settokens[#3]%
\advance\ypos by\height
\putmorphism(\xpos,\ypos)(1,-1)[`\tokend`{#6}]{\height}{\arrowtypec}l%
\puthmorphism(\xpos,\ypos)[\tokena`\tokenb`{#4}]{\height}{\arrowtypea}a%
\advance\xpos by\height
\putvmorphism(\xpos,\ypos)[``{#7}]{\height}{\arrowtyped}m%
\advance\xpos by\height
\putmorphism(\xpos,\ypos)(-1,-1)[``{#8}]{\height}{\arrowtypee}r%
}}

\def\putVtrianglepair{\@ifnextchar <{\putVtrianglepairp}{\putVtrianglepairp%
    <\arrowtypea`\arrowtypeb`\arrowtypec`\arrowtyped`\arrowtypee;\height>}}
\def\Vtrianglepair{\@ifnextchar <{\Vtrianglepairp}{\Vtrianglepairp%
    <\arrowtypea`\arrowtypeb`\arrowtypec`\arrowtyped`\arrowtypee;\height>}}

\def\Vtrianglepairp<#1>[#2;#3`#4`#5`#6`#7]{{%
\settripairparms[#1]%
\settokens[#2]
\xext=\height                  
\width=\height                 
\yext=\height                  
\vertadjust[\tokena`\tokenb`{#4}]
\tempcountd=\tempcounta        
\vertadjust[\tokenb`\tokenc`{#5}]
\ifnum\tempcounta<\tempcountd%
\tempcounta=\tempcountd\fi
\advance \yext by\tempcounta
\botadjust[\tokend``]%
\multiply \xext by2
\leftsladjust[\tokena`\tokend`{#6}]%
\rightsladjust[\tokenc`\tokend`{#7}]%
\begin{picture}(\xext,\yext)(\xoff,\yoff)%
\putVtrianglepairp
<\arrowtypea`\arrowtypeb`\arrowtypec`\arrowtyped`\arrowtypee;\height>%
(0,0)[#2;#3`#4`#5`#6`{#7}]%
\end{picture}%
}}

\def\putCtrianglep<#1>(#2,#3)[#4`#5`#6;#7`#8`#9]{{%
\settriparms[#1]%
\xpos=#2 \ypos=#3
\advance\ypos by\height
\putmorphism(\xpos,\ypos)(1,-1)[``{#9}]{\height}{\arrowtypec}l%
\advance\xpos by\height
\advance\ypos by\height
\putmorphism(\xpos,\ypos)(-1,-1)[#4`#5`{#7}]{\height}{\arrowtypea}l%
{\multiply\height by 2
\putvmorphism(\xpos,\ypos)[`#6`{#8}]{\height}{\arrowtypeb}r}%
}}

\def\putCtriangle{\@ifnextchar <{\putCtrianglep}{\putCtrianglep
    <\arrowtypea`\arrowtypeb`\arrowtypec;\height>}}
\def\Ctriangle{\@ifnextchar <{\Ctrianglep}{\Ctrianglep
    <\arrowtypea`\arrowtypeb`\arrowtypec;\height>}}

\def\Ctrianglep<#1>[#2`#3`#4;#5`#6`#7]{{
\settriparms[#1]
\width=\height                          
\xext=\width                            
\yext=\height                           
\multiply \yext by2 
\topadjust[#2``]
\botadjust[#4``]
\sladjust[#3`#2`{#5}]{\width}
\tempcountd=\tempcounta                 
\sladjust[#3`#4`{#7}]{\width}
\ifnum \tempcounta<\tempcountd          
\tempcounta=\tempcountd\fi              
\advance \xext by\tempcounta            
\advance \xoff by-\tempcounta           %
\rightadjust[#2`#4`{#6}]
\begin{picture}(\xext,\yext)(\xoff,\yoff)%
\putCtrianglep<\arrowtypea`\arrowtypeb`\arrowtypec;\height>%
(0,0)[#2`#3`#4;#5`#6`{#7}]%
\end{picture}%
}}

\def\putDtrianglep<#1>(#2,#3)[#4`#5`#6;#7`#8`#9]{{%
\settriparms[#1]%
\xpos=#2 \ypos=#3
\advance\xpos by\height \advance\ypos by\height
\putmorphism(\xpos,\ypos)(-1,-1)[``{#9}]{\height}{\arrowtypec}r%
\advance\xpos by-\height \advance\ypos by\height
\putmorphism(\xpos,\ypos)(1,-1)[`#5`{#8}]{\height}{\arrowtypeb}r%
{\multiply\height by 2
\putvmorphism(\xpos,\ypos)[#4`#6`{#7}]{\height}{\arrowtypea}l}%
}}

\def\putDtriangle{\@ifnextchar <{\putDtrianglep}{\putDtrianglep
    <\arrowtypea`\arrowtypeb`\arrowtypec;\height>}}
\def\Dtriangle{\@ifnextchar <{\Dtrianglep}{\Dtrianglep
   <\arrowtypea`\arrowtypeb`\arrowtypec;\height>}}

\def\Dtrianglep<#1>[#2`#3`#4;#5`#6`#7]{{
\settriparms[#1]
\width=\height                         
\xext=\height                          
\yext=\height                          
\multiply \yext by2 
\topadjust[#2``]
\botadjust[#4``]
\leftadjust[#2`#4`{#5}]
\sladjust[#3`#2`{#5}]{\height}
\tempcountd=\tempcountd                
\sladjust[#3`#4`{#7}]{\height}
\ifnum \tempcounta<\tempcountd         
\tempcounta=\tempcountd\fi             
\advance \xext by\tempcounta           %
\begin{picture}(\xext,\yext)(\xoff,\yoff)
\putDtrianglep<\arrowtypea`\arrowtypeb`\arrowtypec;\height>%
(0,0)[#2`#3`#4;#5`#6`{#7}]%
\end{picture}%
}}

\def\setrecparms[#1`#2]{\width=#1 \height=#2}%
%

\def\recursep<#1`#2>[#3;#4`#5`#6`#7`#8]{{%
\width=#1 \height=#2
\settokens[#3]
\settowidth{\tempdimen}{$\tokena$}
\ifdim\tempdimen=0pt
  \savebox{\tempboxa}{\hbox{$\tokenb$}}%
  \savebox{\tempboxb}{\hbox{$\tokend$}}%
  \savebox{\tempboxc}{\hbox{$#6$}}%
\else
  \savebox{\tempboxa}{\hbox{$\hbox{$\tokena$}\times\hbox{$\tokenb$}$}}%
  \savebox{\tempboxb}{\hbox{$\hbox{$\tokena$}\times\hbox{$\tokend$}$}}%
  \savebox{\tempboxc}{\hbox{$\hbox{$\tokena$}\times\hbox{$#6$}$}}%
\fi
\ypos=\height
\divide\ypos by 2
\xpos=\ypos
\advance\xpos by \width
\xext=\xpos \yext=\height
\topadjust[#3`\usebox{\tempboxa}`{#4}]%
\botadjust[#5`\usebox{\tempboxb}`{#8}]%
\sladjust[\tokenc`\tokenb`{#5}]{\ypos}%
\tempcountd=\tempcounta
\sladjust[\tokenc`\tokend`{#5}]{\ypos}%
\ifnum \tempcounta<\tempcountd
\tempcounta=\tempcountd\fi
\advance \xext by\tempcounta
\advance \xoff by-\tempcounta
\rightadjust[\usebox{\tempboxa}`\usebox{\tempboxb}`\usebox{\tempboxc}]%
\bfig
\putCtrianglep<-1`1`1;\ypos>(0,0)[`\tokenc`;#5`#6`{#7}]%
\puthmorphism(\ypos,0)[\tokend`\usebox{\tempboxb}`{#8}]{\width}{-1}b%
\puthmorphism(\ypos,\height)[\tokenb`\usebox{\tempboxa}`{#4}]{\width}{-1}a%
\advance\ypos by \width
\putvmorphism(\ypos,\height)[``\usebox{\tempboxc}]{\height}1r%
\efig
}}

\def\recurse{\@ifnextchar <{\recursep}{\recursep<\width`\height>}}

\def\puttwohmorphisms(#1,#2)[#3`#4;#5`#6]#7#8#9{{%
%
\puthmorphism(#1,#2)[#3`#4`]{#7}0a
\ypos=#2
\advance\ypos by 20
\puthmorphism(#1,\ypos)[\phantom{#3}`\phantom{#4}`#5]{#7}{#8}a
\advance\ypos by -40
\puthmorphism(#1,\ypos)[\phantom{#3}`\phantom{#4}`#6]{#7}{#9}b
}}

\def\puttwovmorphisms(#1,#2)[#3`#4;#5`#6]#7#8#9{{%
%
%
%
\putvmorphism(#1,#2)[#3`#4`]{#7}0a
\xpos=#1
\advance\xpos by -20
\putvmorphism(\xpos,#2)[\phantom{#3}`\phantom{#4}`#5]{#7}{#8}l
\advance\xpos by 40
\putvmorphism(\xpos,#2)[\phantom{#3}`\phantom{#4}`#6]{#7}{#9}r
}}

\def\puthcoequalizer(#1)[#2`#3`#4;#5`#6`#7]#8#9{{%
%
\setpos(#1)%
\puttwohmorphisms(\xpos,\ypos)[#2`#3;#5`#6]{#8}11%
\advance\xpos by #8
\puthmorphism(\xpos,\ypos)[\phantom{#3}`#4`#7]{#8}1{#9}
}}

\def\putvcoequalizer(#1)[#2`#3`#4;#5`#6`#7]#8#9{{%
%
%
%
%
\setpos(#1)%
\puttwovmorphisms(\xpos,\ypos)[#2`#3;#5`#6]{#8}11%
\advance\ypos by -#8
\putvmorphism(\xpos,\ypos)[\phantom{#3}`#4`#7]{#8}1{#9}
}}

\def\putthreehmorphisms(#1)[#2`#3;#4`#5`#6]#7(#8)#9{{%
\setpos(#1) \settypes(#8)
\if a#9 %
     \vertsize{\tempcounta}{#5}%
     \vertsize{\tempcountb}{#6}%
     \ifnum \tempcounta<\tempcountb \tempcounta=\tempcountb \fi
\else
     \vertsize{\tempcounta}{#4}%
     \vertsize{\tempcountb}{#5}%
     \ifnum \tempcounta<\tempcountb \tempcounta=\tempcountb \fi
\fi
\advance \tempcounta by 60
\puthmorphism(\xpos,\ypos)[#2`#3`#5]{#7}{\arrowtypeb}{#9}
\advance\ypos by \tempcounta
\puthmorphism(\xpos,\ypos)[\phantom{#2}`\phantom{#3}`#4]{#7}{\arrowtypea}{#9}
\advance\ypos by -\tempcounta \advance\ypos by -\tempcounta
\puthmorphism(\xpos,\ypos)[\phantom{#2}`\phantom{#3}`#6]{#7}{\arrowtypec}{#9}
}}

\def\putarc(#1,#2)[#3`#4`#5]#6#7#8{{%
\xpos #1
\ypos #2
\width #6
\arrowlength #6
\putbox(\xpos,\ypos){#3\vphantom{#4}}%
{\advance \xpos by\arrowlength
\putbox(\xpos,\ypos){\vphantom{#3}#4}}%
\horsize{\tempcounta}{#3}%
\horsize{\tempcountb}{#4}%
\divide \tempcounta by2
\divide \tempcountb by2
\advance \tempcounta by30
\advance \tempcountb by30
\advance \xpos by\tempcounta
\advance \arrowlength by-\tempcounta
\advance \arrowlength by-\tempcountb
\halflength=\arrowlength \divide\halflength by 2
\divide\arrowlength by 5
\put(\xpos,\ypos){\bezier{\arrowlength}(0,0)(50,50)(\halflength,50)}
\ifnum #7=-1 \put(\xpos,\ypos){\vector(-3,-2)0} \fi
\advance\xpos by \halflength
\put(\xpos,\ypos){\xpos=\halflength \advance\xpos by -50
   \bezier{\arrowlength}(0,50)(\xpos,50)(\halflength,0)}
\ifnum #7=1 {\advance \xpos by
   \halflength \put(\xpos,\ypos){\vector(3,-2)0}} \fi
\advance\ypos by 50
\vertsize{\tempcounta}{#5}%
\divide\tempcounta by2
\advance \tempcounta by20
\if a#8 %
   \advance \ypos by\tempcounta
   \putbox(\xpos,\ypos){#5}%
\else
   \advance \ypos by-\tempcounta
   \putbox(\xpos,\ypos){#5}%
\fi
}}

\makeatother

\def\nec{\Box}

\newcommand{\lra}{\longrightarrow}
\newcommand{\ra}{\rightarrow}
\newcommand{\da}{\downarrow}

\def\comp{{\bf Comp}} 
\def\compn{{\bf Comp_n}} 
\def\compnj{{\bf Comp_{n+1}}} 
\def\compnmj{{\bf Comp_{n-1}}} 
\def\compz{{\bf Comp_0}} 
\def\compj{{\bf Comp_1}} 

\def\ccatn{{\bf CCat_n}} 
\def\ccatnj{{\bf CCat_{n+1}}} 
\def\ccatz{{\bf CCat_0}} 
\def\ccatj{{\bf CCat_1}} 

\def\catn{{\bf nCat}} 
\def\catnj{{\bf (n+1)Cat}} 
\def\catnmj{{\bf (n-1)Cat}} 
\def\catz{{\bf 0Cat}} 
\def\catj{{\bf 1Cat}} 

\def\cato{{\bf \o Cat}} 

\def\compnj{{\bf Comp_{n+1}}} 
\def\compnmj{{\bf Comp_{n-1}}} 
\def\compz{{\bf Comp_0}} 
\def\compj{{\bf Comp_1}} 
\def\compt{{\bf Comp_2}} 

\def\comn{{\bf Com_n}} 
\def\comnj{{\bf Com_{n+1}}} 
\def\comz{{\bf Com_0}} 
\def\comj{{\bf Com_1}} 
\def\comt{{\bf Com_2}} 

\def\lk{\langle}
\def\rk{\rangle}

\def\o{\omega}

\pagenumbering{arabic} \setcounter{page}{1}

\title{The category of 3-computads is not cartesian closed}
\author{Mihaly Makkai and Marek Zawadowski\\
Department of Mathematics and Statistics, McGill University,\\
  805 Sherbrooke St., Montr\'eal, PQ, H3A\thinspace2K6, Canada\\[5pt]
Instytut Matematyki, Uniwersytet Warszawski \\
ul. S.Banacha 2, 00-913 Warszawa, Poland}
\date{February 26, 2008}

\maketitle

\begin{abstract} We show, using \cite{CarbJohn} and Eckmann-Hilton
argument, that the category of $3$-computads is not cartesian
closed. As a corollary we get that neither the category of all
computads nor the category of $n$-computads, for $n>2$, do form
locally cartesian closed categories, and hence elementary toposes.
\end{abstract}
\section{Introduction}

S.H. Schanuel (unpublished) made an observation, c.f.
\cite{CarbJohn}, that the category of $2$-computads $\comp_{\bf
2}$ is a presheaf category. We show below that neither the
category of computads nor the categories $n$-computads, for $n>2$,
are locally cartesian closed. This is in contrast with a remark in
\cite{CarbJohn} on page 453, and an explicit statement in
\cite{Batanin} claiming that these categories are presheaves
categories. Note that some interesting subcategories of computads,
like many-to-one computads, do form presheaf categories, c.f. \cite{HMP}, \cite{HMZ}.

We thank the anonymous referee for comments that helped to clarify
the exposition of the example. The diagrams for this paper were prepared
with a help of {\em catmac} of Michael Barr.

\section{Computads}
Computads were introduced by R.Street in \cite{Street}, see also \cite{Batanin}. Recall
that a computad is an $\o$-category that is levelwise free.
Below we recall one of the definitions.

Let $\catn$ be the category of $n$-categories and $n$-functors between them,
$\cato$ be the category of $\o$-categories and $\o$-functors between them.
We have the obvious truncation functors \[ tr_{n-1} : \catn\lra \catnmj \]
By $\compn$ we denote the category of $n$-computads, a non-full subcategory
of the category $\catn$. By $\ccatn$ we denote the non-full subcategory of $\catn$,
whose objects are 'computads up to the level $n-1$', i.e.  an $n$-functor $f:A\ra B$ is a morphism
in $\ccatn$ if and only if $tr_{n-1}(f):tr_{n-1}(A)\ra tr_{n-1}(B)$ is a morphism in $\compnmj$.
Clearly $\ccatn$ is defined as soon as $\compnmj$ is defined. The categories
$\compn$ and $n$-comma category $\comn$ are defined below.

The categories $\compz$, $\ccatz$ and $\comz$ are equal to $Set$, the category of sets.
We have an adjunction
\begin{center}
\xext=900 \yext=200
\begin{picture}(\xext,\yext)(\xoff,\yoff)
\putmorphism(0,100)(1,0)[\comz`\ccatz`]{900}{0}a
\putmorphism(0,50)(1,0)[\phantom{\comz}`\phantom{\ccatz}`U_0]{900}{-1}b
\putmorphism(0,150)(1,0)[\phantom{\comz}`\phantom{\ccatz}`F_0]{900}{1}a
\end{picture}
\end{center}
with both functors being the identity on $Set$, $F_0\dashv U_0$. $\compz$ is
the image of $\comz$ under $F_0$.

$\comj$ is the category of graphs, i.e. an object of $\comj$ is a pair of sets
and a pair of functions between them $\lk d,c:E\ra V\rk$.
$\ccatj$ is simply $\bf Cat$, the category of all small categories.
The forgetful functor $U_1$ (forgetting compositions and identities) has a left
adjoint $F_1$ 'the free category (over a graph)' functor
\begin{center}
\xext=900 \yext=200
\begin{picture}(\xext,\yext)(\xoff,\yoff)
\putmorphism(0,100)(1,0)[\comj`\ccatj`]{900}{0}a
\putmorphism(0,50)(1,0)[\phantom{\comj}`\phantom{\ccatj}`U_1]{900}{-1}b
\putmorphism(0,150)(1,0)[\phantom{\comj}`\phantom{\ccatj}`F_1]{900}{1}a
\end{picture}
\end{center}
We have a diagram
\begin{center}
\xext=1200 \yext=750
\begin{picture}(\xext,\yext)(\xoff,\yoff)
 \setsqparms[1`1`1`1;1200`600]
 \putsquare(0,50)[\comj`\ccatj`\comz`\ccatz;F_1`tr_0'`tr_0`F_0]
 \put(420,300){\mbox{$\compz$}}
 \put(100,100){\vector(2,1){370}}
  \put(200,220){\mbox{$F_0$}}
 \put(700,270){\vector(2,-1){350}}
  \put(850,220){\mbox{$\iota_0$}}
 \put(100,580){\vector(2,-1){370}}
 \put(230,520){\mbox{$tr_0$}}
  \put(1070,575){\vector(-2,-1){370}}
 \put(800,520){\mbox{$tr_0$}}
\end{picture}
\end{center}
where three triangles commute, moreover the left triangle and the outer square
commute up to an isomorphism.  $tr_1$ and $tr'_1$ are the obvious truncation morphisms.
Then we define the category of $1$-computads $\compj$ as the essential (non-full) image
of the functor $F_1$ in $\ccatj$, i.e. $1$-computads are the free categories over graphs
and computad maps between them are functors sending indets (=indeterminates=generators) to indets.

Now suppose that we have an adjunction $U_n\dashv F_n$
\begin{center}
\xext=900 \yext=700
\begin{picture}(\xext,\yext)(\xoff,\yoff)
\putmorphism(0,100)(1,0)[\comn`\ccatn`]{1200}{0}a
\putmorphism(0,50)(1,0)[\phantom{\comn}`\phantom{\ccatn}`U_n]{1200}{-1}b
\putmorphism(0,150)(1,0)[\phantom{\comn}`\phantom{\ccatn}`F_n]{1200}{1}a

 \put(420,400){\mbox{$\compn$}}
 \put(100,200){\vector(2,1){370}}
  \put(200,320){\mbox{$F_n$}}
 \put(700,370){\vector(2,-1){350}}
  \put(850,320){\mbox{$\iota_n$}}
\end{picture}
\end{center}
and $\compn$ is defined as the the essential (non-full) image
of the functor $F_n$ in $\ccatn$. We define the $n$-parallel pair functor
\begin{center}
\xext=800 \yext=150
\begin{picture}(\xext,\yext)(\xoff,\yoff)
\putmorphism(0,50)(1,0)[\Pi_n:\compn`Set`]{800}{1}b
\end{picture}
\end{center}
such that
\[ \Pi_n(A)=\{ \lk a,b\rk |\; a,b\in A_n,\; d(a)=d(b),\; c(a)=c(b) \}  \]
for any $n$-computad $A$. The $(n+1)$-comma category $\comnj$ is the category $Set\da \Pi_n$.
Thus an object in $\comnj$ is a pair $( A, \lk d,c  \rk : X \ra \Pi_n(A)$, such that $A$ is
an $n$-computad $X$ is a set of $(n+1)$-indets and $\lk d,c  \rk$ is a function associating $n$-domains
and $n$-codomains. The forgetful functor $U_{n+1}:\ccatnj\lra \comnj$ (forgetting compositions
and identities at the level $n+1$) creates limits and satisfies the solution set condition. Thus it has a left adjoint $F_{n+1}$.
We get a diagram
\begin{center}
\xext=1200 \yext=750
\begin{picture}(\xext,\yext)(\xoff,\yoff)
 \setsqparms[1`1`1`1;1200`600]
 \putsquare(0,50)[\comnj`\ccatnj`\comn`\ccatn;F_{n+1}`tr_n'`tr_n`F_n]
 \put(420,300){\mbox{$\compn$}}
 \put(100,100){\vector(2,1){370}}
  \put(190,220){\mbox{$F_n$}}
 \put(700,270){\vector(2,-1){350}}
  \put(850,220){\mbox{$\iota_n$}}
 \put(100,580){\vector(2,-1){370}}
 \put(250,520){\mbox{$tr_{n}$}}
 \put(1070,575){\vector(-2,-1){370}}
 \put(760,520){\mbox{$tr_{n}$}}
\end{picture}
\end{center}
where three triangles commute, moreover the left triangle and the outer square
commute up to an isomorphism.  $tr_{n}$ are the obvious truncation functors and $tr'_{n}$
is a truncation functor that at the level $n$ leaves the indets only.
Then we define the category of $(n+1)$-computads $\compnj$ as the essential (non-full) image
of the functor $F_{n+1}$ in $\ccatnj$, i.e. $(n+1)$-computads are the free $(n+1)$-categories
over $(n+1)$-comma categories and $(n+1)$-computad maps between them are $(n+1)$-functors sending
indets to indets. The category of computads $\comp$ is a (non-full) subcategory of the category of
$\o$-categories and $\o$-functors $\cato$ such, that for each $n$, the truncation of objects and
morphisms to $\catn$ is in $\compn$. As $F_n: \comn \ra \ccatn$ is faithful and full on isomorphisms,
after restricting the codomain we get an equivalence of categories $F_n: \comn \ra \compn$.

{\em Notation.} If $A$ is a computad then $A_n$ denotes the set of $n$-cells of $A$ and $|A|_n$ denotes
the set of $n$-indets of $A$.

The truncation functor $tr_n : \compnj\lra\compn$ has both adjoints \mbox{$i_n\dashv tr_n \dashv f_n$}
\begin{center}
\xext=900 \yext=250
\begin{picture}(\xext,\yext)(\xoff,\yoff)
\putmorphism(0,100)(1,0)[\compnj`\compn`tr_n]{900}{1}a
\putmorphism(0,50)(1,0)[\phantom{\compnj}`\phantom{\compn}`i_n]{900}{-1}b
\putmorphism(0,220)(1,0)[\phantom{\compnj}`\phantom{\compn}`f_n]{900}{-1}a
\end{picture}
\end{center}
where
\[ i_n(A)= F_{n+1}(A, \emptyset \ra \Pi_n(A)) \]
and
\[ f_n(A)= F_{n+1}(A, id_{\Pi_n(A)}: \Pi_n(A) \ra \Pi_n(A)) \]
for $A$ in $\compn$. This shows that $tr_n$ preserves limits and colimits. The colimits in $\compnj$
are calculated in $\catnj$ but the limits in $\compnj$ are more involved. It is more convenient to
describe them in $\comnj$ and then apply the functor $F_{n+1}$. If $H:{\cal J} \ra \comnj$ is a functor and
$P$ is the limit of its truncation $tr_n\circ H$ to $\compn$ then $Lim\, H$, the limit of $H$,
truncated to $\compn$ is $P$ and the $(n+1)$-indets $|Lim\, H|_{n+1}$ of $Lim\, H$ are as follows
\[  |Lim\, H|_{n+1}=\{ \lk a_i\rk_{i\in {\cal J} } |\; a_i\in |H(i)|_{n+1},\;
\lk d(a_i)\rk_{i\in {\cal J} },\, \lk c(a_i)\rk_{i\in {\cal J} } \in P_n  \}  \]

The terminal object $1_n$ in $\compn$ is quite complicated, for $n\geq 2$. However the $\comt$ part of $1_2$
is still easy to describe. $1_2$ has one $0$-indet $x$ and one $1$-indet $\xi :x\ra x$. Thus the $1$-cells
can be identified with finite (possibly empty) strings of of arrows:
\begin{center}
\xext=1700 \yext=150
\begin{picture}(\xext,\yext)(\xoff,\yoff)
\putmorphism(500,50)(1,0)[x`x`\xi]{300}{1}a
\putmorphism(800,50)(1,0)[\phantom{x}`x`\xi]{300}{1}a
\put(1200,40){\mbox{$\ldots$}}
\putmorphism(1400,50)(1,0)[x`x`\xi]{300}{1}a
\put(0,30){\mbox{$x,$}}
\end{picture}
\end{center}
or simply with elements of $\o$.  The set $|1_2|_{2}$ of $2$-indets in $1_2$ contains exactly one indet for
every pair of strings. The first element of such a pair is the domain of the indet and the second element
of the pair is the codomain of the indet. Thus $|1_2|_{2}$ can be identified with the set $\o\times\o$.
In particular $\lk 0,0\rk$ correspond to the only indet from $id_x$ to $id_x$ ($id_x$ is the identity on $x$).
The description of all $2$-cells in $1_2$ is more involved but we don't need it here.

\section{The counterexample}

\begin{lemma}\label{notcc}
$\bf \comp_3$ is not cartesian closed.
\end{lemma}

{\it Proof.}~
As it was noted in Lemma 4.2  \cite{CarbJohn}, the functor $\Pi_2$
factorizes as
\begin{center} \xext=1700 \yext=100
\begin{picture}(\xext,\yext)(\xoff,\yoff)
 \putmorphism(0,0)(1,0)[{\compt}`Set\da \Pi_2(1_2)`\widehat{\Pi_2}]{850}{1}a
 \putmorphism(850,0)(1,0)[\phantom{Set\da \Pi_2(1_2)}`Set`\Sigma]{850}{1}a
\end{picture}
\end{center}
where $\widehat{\Pi_2}(A)=\Pi_2(!:A\ra 1_2)$, and $\Sigma(b:B\ra\Pi_2(1_2))=B$,
for $A$ in $\compt$ and $b$ in $Set\da \Pi_2(1_2)$. Moreover, the category $Set\da \Pi_2$,
which is equivalent to $\bf \comp_3$, is also equivalent to $(Set\da \Pi_2(1_2))\da \widehat{\Pi_2}$.
Now, as $\compt$ and $Set\da \Pi_2(1_2)$ are cartesian closed categories with initial objects
(in fact both categories are presheaf toposes) and $\widehat{\Pi_2}$ preserves the terminal object,
by Theorem 4.1 of \cite{CarbJohn}, $\bf Comp_3$ is a cartesian closed category if and only if  $\widehat{\Pi_2}$
preserves binary products. We finish the proof by showing that $\widehat{\Pi_2}$
does not preserves the binary products.

Let $A$ be a $2$-computad with one $0$-cell $x$, one $1$-cell
$id_x$ the identity on $x$ (no $1$-indets).  Moreover $A$ has as $2$-cells all
cells generated by the two indeterminate $2$-cells $a_1,a_2:id_x\ra
id_x$. Thus, by Eckmann-Hilton argument, any $2$-cell in $A$ is of
form $a_1^m\circ a_2^n$, for $m,n\in\omega$ (if $m=n=0$ then
$a_1^m\circ a_2^n=id_{id_x}$). Let $B$ be a $2$-computad
isomorphic to $A$ with indeterminate $2$-cells  $b_1$, $b_2$. Let
$x$ be the unique $0$-cell in $1_2$, $c$ be the only indeterminate
$2$-cell in $1_2$ that has $id_x$ as its domain and codomain and
$C$ a subcomputad of $1_2$ generated by $c$. The unique maps of
$2$-computads $! :A\ra 1_2$ and $! :B\ra 1_2$ sends $a_i$ and
$b_i$ to $c$, for $i=1,2$. Thus they factor through $C$ as $\alpha
:A\ra C$ and $\beta :B\ra C$, respectively. The $2$-computad $C$
does not play a crucial role in the counterexample but it makes
the explanations simpler.

Let us describe the product $A\times B$ in $\compt$. The $0$-cell
and $1$-cells are as in $A$, $B$ and $C$. As there is only one $1$-cell
$id_x$ in $A\times B$, the compatibility condition for domain and codomains
of $2$-indets is trivially satisfied, and   the set $2$-indets of $A\times B$
is just the product of $2$-indets of $A$ and $B$, i.e.
\[ |A\times B|_2= \{ \lk a_i,b_j\rk |\; i,j=1,2\} \]
and the set of all $2$-cells of $A\times B$ is
\[ (A\times B)_2= \{ \lk a_1,b_1\rk^{n_1}\circ\lk a_1,b_2\rk^{n_2}\circ\lk a_2,b_1\rk^{n_3}
\circ\lk a_2,b_2\rk^{n_2} |\; n_1,n_2,n_3,n_4\in \o \} \]
The projections
\begin{center}
\xext=1200 \yext=120
\begin{picture}(\xext,\yext)(\xoff,\yoff)
\putmorphism(0,50)(1,0)[A`A\times B`\pi_1]{600}{-1}a
\putmorphism(600,50)(1,0)[\phantom{A\times B}`B`\pi_B]{600}{1}a
\end{picture}
\end{center}
are defined as the only $2$-functors such that  $\pi_A(a_i,b_j)=a_i$ and $\pi_A(a_i,b_j)=b_j$,
for $i,j=1,2$. Thus we have a commuting square
\begin{center}
\xext=600 \yext=900
\begin{picture}(\xext,\yext)(\xoff,\yoff)
\settriparms[1`1`0;300]
 \putAtriangle(0,550)[A\times B`A`B;\pi_A`\pi_B`]
 \settriparms[0`1`1;300]
 \putVtriangle(0,250)[\phantom{A}`\phantom{B}`C;`\alpha`\beta]
\putmorphism(300,260)(0,-1)[`1_2`m]{280}{2}r
\put(0,500){\line(0,-1){380}}
\put(600,500){\line(0,-1){380}}
\put(0,120){\vector(2,-1){240}}
\put(600,120){\vector(-2,-1){240}}
\put(-40,200){\mbox{$!$}}
\put(620,200){\mbox{$!$}}
\put(1920,400){\mbox{$(*)$}}
\end{picture}
\end{center}
As $C$ is a subobject of the terminal object $A\times B$ is $A\times_C B$ and $A\times_{1_2} B$,
i.e. both inner and outer squares in the above diagram are pullbacks.

Since all the $2$-cells in $A$, $B$, $C$ and $A\times B$ are parallel we have
\[ \Pi_2(A)=A_2\times A_2,\;\;\;\; \Pi_2(B)=B_2\times B_2,\;\;\;\; \Pi_2(C)=C_2\times C_2,\]
and\[ \Pi_2(A\times B)=(A\times B)_2\times (A\times B)_2  \]
$\widehat{\Pi_2}$ preserves the product of $A$ and $B$ if in the diagram $(**)$ below,
which is the application of $\Pi_2$ to the diagram $(*)$ above, the outer square is
a pullback in $Set$
\begin{center}
\xext=600 \yext=1350
\begin{picture}(\xext,\yext)(\xoff,\yoff)
\settriparms[1`1`0;450]
 \putAtriangle(150,850)[(A\times B)_2\times(A\times B)_2`A_2\times A_2`B_2\times B_2;\Pi_2(\pi_A)`\Pi_2(\pi_B)`]
 \settriparms[0`1`1;450]
 \putVtriangle(150,400)[\phantom{A_2\times A_2}`\phantom{B_2\times B_2}`C_2\times C_2;`\Pi_2(\alpha)`\Pi_2(\beta)]

\putmorphism(600,400)(0,-1)[\phantom{C_2\times C_2}`\Pi_2(1_2)`\Pi_2(m)]{400}{1}r
\put(0,750){\line(0,-1){500}}
\put(1200,750){\line(0,-1){500}}
\put(0,250){\vector(2,-1){400}}
\put(1200,250){\vector(-2,-1){400}}
\put(-240,350){\mbox{$\Pi_2(!)$}}
\put(1220,350){\mbox{$\Pi_2(!)$}}
\put(1850,400){\mbox{$(**)$}}
\end{picture}
\end{center}
As $\Pi_2(m)$ is mono, the outer square in $(**)$ is a pullback in $Set$ if and only
if the inner square in $(**)$ is a pullback in $Set$. We have \[ \Pi_2(\pi_A)=(\pi_A)_2\times(\pi_A)_2,
\;\;\; \Pi_2(\pi_B)=(\pi_B)_2\times(\pi_B)_2,\]  \[ \Pi_2(\alpha)=\alpha_2\times\alpha_2,
\;\;\; \mbox{and}\;\;\; \Pi_2(\beta)=\beta_2\times\beta_2.\] Hence the inner square in
$(**)$ is a pullback if and only
if the square $(***)$ below
\begin{center}
\xext=600 \yext=650
\begin{picture}(\xext,\yext)(\xoff,\yoff)
\settriparms[1`1`0;300]
 \putAtriangle(0,300)[(A\times B)_2`A_2`B_2;(\pi_A)_2`(\pi_B)_2`]
 \settriparms[0`1`1;300]
 \putVtriangle(0,0)[\phantom{(A)_2}`\phantom{(B)_2}`(C)_2;`\alpha_2`\beta_2]
\put(1820,400){\mbox{$(***)$}}
\end{picture}
\end{center}
is a pullback. But $(***)$ is not a pullback in $Set$. The two $2$-cells
\[ \lk a_1,b_1\rk\circ \lk a_2,b_2\rk, \;\; \mbox{and} \;\;\lk a_1,b_2\rk\circ
\lk a_2,b_1\rk\]
 in $A\times B$ are different since they are compositions of different indets.
On the other hand
\[ (\pi_A )_2((a_1,b_1)\circ
(a_2,b_2))=a_1\circ a_2 = (\pi_A )_2((a_1,b_2)\circ (a_2,b_1))\]
and
\[ (\pi_B )_2((a_1,b_1)\circ
(a_2,b_2))=b_1\circ b_2=b_2\circ b_1= (\pi_B )_2((a_1,b_2)\circ (a_2,b_1))\]
i.e. they agree on both projections and hence $(***)$ is not a pullback.
Thus $\widehat{\Pi_2}$ does not preserve binary products, as required. $~~\Box$

\begin{theorem}\label{thm}
 The category of computads $\comp$ and the categories
of $n$-computads $\compn$, for $n>2$, are not locally cartesian
closed.
\end{theorem}

{\it Proof.}~ The slice categories $\comp\da 1_3$, as well as
$\compn\da 1_3$, for $n>2$, are equivalent to $\comp_{\bf 3}$,
where $1_3$ is the terminal object in $\comp_{\bf 3}$ lifted (by
adding suitable identities) to the category of appropriate
computads. As, by Lemma \ref{notcc}, $\compn\da 1_3$ is not
cartesian closed we get the theorem. $~\Box$

{\em Remark.} In particular the categories mentioned in the above
theorem are not presheaf (or even elementary) toposes.


\begin{thebibliography}{CWMW}
\bibitem[B]{Batanin}
 \frenchspacing
M.Batanin,  {\em Computads for finitary monads on globular sets}.
Contemporary Mathematics, vol 230, (1998), pp. 37-57.

\bibitem[CJ]{CarbJohn} A. Carboni, P.T. Johnstone {\em Connected limits,
familial representability and Artin glueing}. Math. Struct. in
Comp Science, vol 5. (1995), pp. 441-459.

\bibitem[HMZ]{HMZ} V. Harnik, M. Makkai, M. Zawadowski,
{\em Multitopic sets are the same as many-to-one computads}. Preprint 2002.

\bibitem[HMP]{HMP} C. Hermida, M. Makkai, J. Power, {\em On weak
higher dimensional categories, I} Parts 1,2,3, J. Pure and Applied
Alg.  153 (2000), pp. 221-246, 157 (2001), pp. 247-277, 166
(2002), pp. 83-104.

\bibitem[S]{Street}
R.Street,  {\em Limits indexed by category valued 2-functors}. J.
Pure Apll. Alg, vol 8, (1976), pp. 149-181.
\end{thebibliography}
\end{document}